\documentclass[11pt]{article}
\usepackage[papersize={8.5in,11in},top=0.9in,left=1in,right=1in,bottom=0.92in]{geometry}
\usepackage{cite}

\usepackage{url,amssymb,graphics,amsmath,amsthm,amsopn,amstext,amsfonts}
\usepackage{hyperref}
\usepackage{graphicx}
\usepackage{caption}
\usepackage{wrapfig}
\usepackage{times}
\usepackage[noend]{algpseudocode}
\usepackage{algorithm}

\newtheorem{lemma}{Lemma}
\newtheorem{proposition}{Proposition}
\newtheorem{theorem}{Theorem}
\newtheorem{corollary}{Corollary}
\newtheorem{remark}{Remark}

\begin{document}


\title{Knot Locating in Piecewise Linear Approximation}
\author{Carlos Ugaz\footnote{Email: cugaz@precima.com.}, Lanshan Han\footnote{Corresponding author, Email: lhan2@precima.com.}, and Alvin Lim \footnote{Email: alim@precima.com.} \footnote{The authors are with Precima, Inc., 8600 W Bryn Mawr Ave., Suite 1000N, Chicago, IL 60631.}}

\maketitle
\abstract{Many separable nonlinear optimization problems can be approximated by their nonlinear objective functions with piecewise linear functions. A natural question arising from applying this approach is how to break the interval of interest into subintervals (pieces) to achieve a good approximation. We present formulations to optimize the location of the knots. We apply a sequential quadratic programming method and a spectral projected gradient method to solve the problem. We report numerical experiments to show the effectiveness of the proposed approaches. }

\section{Background and Motivations}
In practice, we often need to solve separable nonlinear optimization problems of the following format:
\begin{equation}\label{eq:separable_nonlinear}
\begin{array}{rl}
\min & \displaystyle{\sum_{i=1}^n f_i(y_i) } \\[5pt]
\mbox{s.t.} & A\boldsymbol y \, \leq \, \boldsymbol b, \\
& \boldsymbol l \, \leq \, \boldsymbol y \, \leq \, \boldsymbol u,
\end{array}
\end{equation}
where $\boldsymbol y =[y_1,\cdots,y_n]^T \in \mathbb R^n$ is the vector of decision variables, $\boldsymbol l \in \mathbb R^n$ and $\boldsymbol u \in \mathbb R^n$ are vectors of the lower and upper bounds, $A\in \mathbb R^{m\times n}$ is the constraint matrix, and $\boldsymbol b\in \mathbb R^{m}$ is the right-hand side. Some prominent examples of the above optimization problem are utility maximization, marketing mix optimization, and others. An example of utility maximization is provided in the next subsection. An important feature of (\ref{eq:separable_nonlinear}) is that the objective function is the sum of a set of univariate functions, and hence, is separable.

To solve this type of separable nonlinear optimization problem, one approach is to first approximate the function $f_i(\cdot)$'s by piecewise linear (PL) functions and then solve the problem as a linear program (when all the $f_i(\cdot)$'s are convex) or a mixed integer linear program (when $f_i(\cdot)$'s are not necessarily convex) \cite{JBard03}. The motivation of utilizing the piecewise linear approximation is the fact that linear programming or mixed integer linear programming solvers are arguably more mature and more accessible to industrial users compared to their nonlinear counterparts. Therefore, while we lose some accuracy doing the problem approximation, we can take advantage of more mature computational tools to solve very large problems, which occur very often in practice, in reasonable computing time. Moreover, in practice, the functions $f_i(\cdot)$'s are often obtained by performing a nonlinear regression using collected data, and hence, is subject to intrinsic inaccuracy. Therefore, it is not necessary to pursue extreme accuracy in the optimization. It is usually acceptable to use PL approximations.

Another reason for the piecewise linear approximation of the functions $f_i(\cdot)$'s is the faster speed of obtaining an approximate solution to (\ref{eq:separable_nonlinear}) via linear programming versus nonlinear programming. While it takes some time to obtain the piecewise approximation of the $f_i(\cdot)$'s, the approximation process is done infrequently and should coincide with the update of the $f_i(\cdot)$'s via nonlinear regression. In contrast, the optimization of (\ref{eq:separable_nonlinear}) with the piecewise approximation of the $f_i(\cdot)$'s is required repeatedly for various business scenarios and at various times until the next update of the $f_i(\cdot)$'s. Therefore, it makes sense to invest time to obtain the piecewise approximation for speedier solutions of (\ref{eq:separable_nonlinear}).

Suppose we are interested in approximating a univariate function $f(x)$ within range $[a,b]$ using a piecewise linear function. To construct the PL approximation, we are given a set of $n$ break points (knots), $x_1,\cdots, x_n$ satisfying $a \leq x_1 \leq x_2 \cdots \leq x_n \leq b$. To simplify the expression, we let $x_0 = a$ and $x_{n+1}=b$ and we use the following PL function to approximate $f(x)$:
\begin{equation}\label{eq:PL_format}
\widehat f(x) \, = \, \alpha_i x + \beta_i, \,\, \mbox{if } x \in [x_i,x_{i+1}],
\end{equation}
where
\begin{equation}\label{eq:alpha_beta}
\begin{array}{rcl}
\alpha_i & = & \displaystyle{\frac{f(x_{i+1}) -f(x_i)}{x_{i+1}-{x_i}}},\,\,\, i=0,\cdots,n; \\[5pt]
\beta_i & = & \displaystyle{\frac{x_{i+1}f(x_{i}) - x_if(x_{i+1})}{x_{i+1}-x_i}},\, \,\, i=0,\cdots,n.
\end{array}
\end{equation}

As we can see, once the knots are given, the approximation is determined, and so is the approximation error. It is obvious that if we increase the number of knots, we can certainly refine the approximation. However, increasing the number of knots leads to larger optimization problems and longer computing times. For a practical problem involving possibly millions of $f_i(\cdot)$'s, this increase in runtime can be dramatic. Therefore, in this paper we aim to refine the approximation by properly locating the knots without increasing their number.
\begin{center}
\begin{minipage}[t]{0.4\textwidth}
\includegraphics[width=\textwidth]{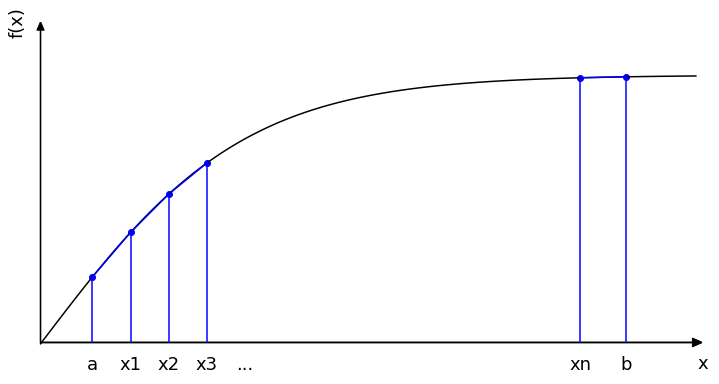}
\end{minipage}
\begin{minipage}[t]{0.4\textwidth}
\includegraphics[width=\textwidth]{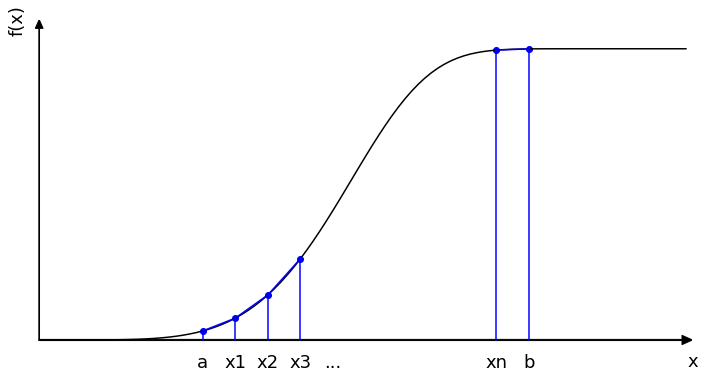}
\end{minipage}\\
Figure 1: Piecewise linear approximation of a function. Left: C-shaped curve. Right: S-shaped curve.
\label{fig:approximation}
\end{center}

The problem of optimizing piecewise linear approximations is not very well studied in the literature. In the case where we are given candidate knots to choose from, the problem can be converted to a network flow problem \cite{ahuja1993network}. However, to the best of our knowledge, the case where the candidate knots is not given has not been studied. This paper is dedicated to applying an optimization approach in locating the knots so that the approximation error, measured by certain norm of the difference between the original function and its approximation, is minimized. The rest of the paper is organized as follows: in the next subsection, we provide an example from internet service management to further motivate applying PL approximation in solving separable nonlinear programs; we present the assumptions, formulations, as well as some mathematical properties in Section~\ref{sec:pre}; a gradient projection algorithm is presented in Section~\ref{sec:alg} followed by experiment results in Section~\ref{sec:exp}; and we conclude the paper in Section~\ref{sec:con} with some closing remarks.

\subsection{A Rate Control Problem in Multi-Class Internet Service}
We consider a service rate control problem presented in \cite{lee2005non}. The purpose of rate control in internet service is to provide satisfactory service to the users while at the same time alleviate network congestion. Typically, an internet service provider can adjust its data transmission rates according to congestion levels within its network. Hence, by allocating appropriate rates for different users, one can maximize user satisfaction while keeping network congestion levels low. A utility function is widely used to measure user satisfaction or Quality of Service (QoS).

For different types of services on the internet, the shapes of the utility functions are different. Roughly speaking, the utility functions can be divided into two classes. The first class pertains to traditional data services, such as file transfer and email. The QoS of this class of services gracefully degrades as the data transmission rate decreases when network congestion is present. The utility functions for this class of services is hence concave (C-chaped) as shown on the left in Figure \ref{fig:approximation}. The second class pertains to streaming video and audio service, in which QoS drops dramatically when the data transmission rate is below a certain threshold. The utility functions for this second class of services can often be assumed to be sigmoidal as shown on the right in Figure \ref{fig:approximation}.

Consider a network with $n$ users indexed by $i$ and $k$ links indexed by $l$. Each link $l$ has capacity $c_l$ and each user has a utility function of his/her data transmission rate, denoted by $U_i(x_i)$, where $x_i$ is user $i$'s data transmission rate. Let $m_i$ be the maximum data transmission rate user $i$ can receive. We define the user-link incidence matrix $A \in \mathbb R^{k\times n}$ with the elements
$$
a_{li} \, = \, \left\{
\begin{array}{ll}
1, & \mbox{if user }i \mbox{ is connected to link }l\\
0, & \mbox{otherwise}
\end{array}
\right..
$$
Letting $\boldsymbol x = \left[x_1,\cdots,x_n\right]^T$ be the vector of the data transmission rate, $\boldsymbol c = \left[c_1,\cdots,c_k\right]^T$ be the vector of link capacities, and $\boldsymbol m=\left[m_1,\cdots,m_n\right]^T$ be the vector of user maximum transmission rates, the rate control problem can be formulated as the following optimization problem:
\begin{equation}\label{eq:Rate_control}
\begin{array}{rl}
\displaystyle{\max_{\boldsymbol x}} & \displaystyle{\sum_{i=1}^N U_i(x_i) } \\[5pt]
\mbox{subject to} & A \boldsymbol x \, \leq \, \boldsymbol c, \\[5pt]
& 0 \, \leq \, \boldsymbol x \, \leq \, \boldsymbol m.
\end{array}
\end{equation}
This is clearly a separable optimization of potentially very large size (depending on the number of users). One may apply the PL approximation approach to find an approximate solution in reasonable time.
\section{Mathematical Models and Properties}\label{sec:pre}
In this section, we present our formulations to study the knot locating problem. The mathematical properties of the formulations are also studied. When the curve is concave(convex), the PL approximation is always an under(over) approximation, leading to a simplified formulation and analysis. In Subsection~\ref{subsec:concave}, we mainly focus on the concave case while the convex case can be studied similarly. In Subsection~\ref{subsec:general} we move on to general curves without concavity/convexity assumptions. Throughout the paper, we use the following convention in our notations: lower case letters represent scalars, bold lower case letters represent vectors, and upper case letters represent matrices.
\subsection{Approximating Concave Curves}\label{subsec:concave}
We study concave curves (also referred to as C-shaped curves in many applications) in this subsection. Mathematically, given a curve $f(x)$ and a range of interest $[a,b]$ with $b>a$, we assume that there exists an $\epsilon>0$, such that
\begin{itemize}
\item[(a)] $f(x)$ is twice differentiable on $(a-\epsilon,b+\epsilon)$. Moreover, $f(x)>0$, $f''(x)<0$, for all $x\in (a-\epsilon,b+\epsilon)$;
\item[(b)] $f(x)$ is an increasing curve, i.e. $f'(x)>0$ for all $x\in (a-\epsilon,b+\epsilon)$;
\item[(c)] $f(x)$ is concave, i.e., $f''(x)<0$ for all $x\in (a-\epsilon,b+\epsilon)$.
\end{itemize}
Note that we can assume without loss of generality that $f(x)$ is positive in the interval of interest. In fact, if this condition is not satisfied, we can simply add a constant to $f(x)$ to satisfy the positivity condition without changing the nature of the problem we are considering in this paper. Let $x_1, x_2, \cdots, x_n$ be the knots in $[a,b]$. To simplify the notation, also let $x_0 = a$ and $x_{n+1} = b$. The area under the curve can be computed by
$$A \, = \, \int_a^b f(x) dx.$$ The area under the PL approximation is given by
$$\widetilde A \, =\, \frac{1}{2}\sum_{i=0}^n (x_{i+1} - x_i) (f(x_{i+1}) + f(x_i)).$$ Since the PL approximation is a piecewise under approximation, as illustrated on the left in Figure~\ref{fig:approximation}, we can measure the error by
\begin{equation}\label{eq:concave_measure}
A - \widetilde A \, =  \, \int_a^b f(x) dx - \frac{1}{2}\sum_{i=0}^n (x_{i+1} - x_i) (f(x_{i+1}) + f(x_i)).
\end{equation}
Therefore, given the number of knots, to find the best approximation, we need to solve the following optimization problem:
\begin{equation}\label{eq:error_minimization}
\begin{array}{rl}
\displaystyle{ \min_{x_1,\cdots, x_n} } & \displaystyle{ -\frac{1}{2}\sum_{i=0}^n (x_{i+1} - x_i) (f(x_{i+1}) + f(x_i)) }\\[10pt]
\mbox{s.t.} & a \, = \, x_0 \, \leq \, x_1 \, \leq \, \cdots \,\leq \, x_n \, \leq x_{n+1} \, =\, b
\end{array}
\end{equation}
Let $\boldsymbol x = [x_1,\cdots,x_n]^T$. To analyze the properties of this minimization problem (\ref{eq:error_minimization}), we write the constraints in a different form and introduce a multiplier $\lambda_i$ for each constraint as follows:
\begin{equation}
\begin{array}{rcc}
\displaystyle{ \min_{x_1,\cdots, x_n} } & \phi( \boldsymbol x) \, \triangleq \, \displaystyle{ -\frac{1}{2}\sum_{i=0}^n (x_{i+1} - x_i) (f(x_{i+1}) + f(x_i)) }\\[10pt]
\mbox{s.t.} & x_0 - x_1 \, \leq \, 0 & (\lambda_0) \\
& x_1 - x_2 \, \leq \, 0 & (\lambda_1)\\
& \vdots & \vdots\\
& x_n - x_{n+1} \, \leq \, 0 & (\lambda_n)
\end{array}
\end{equation}
Let $\boldsymbol \lambda = [\lambda_0,\cdots,\lambda_n]^T$. Define the Lagrangian function as follows:
\begin{equation}
L(\boldsymbol x,\boldsymbol \lambda) \, \triangleq \, -\frac{1}{2}\sum_{i=0}^n (x_{i+1} - x_i) (f(x_{i+1}) + f(x_i)) + \sum_{i=0}^n \lambda_i(x_i-x_{i+1}).
\end{equation}
Since the constraints are linear, the constraint qualification holds. Therefore, if $\boldsymbol x$ is a (local) optimal solution of (\ref{eq:error_minimization}), then there must exist $\boldsymbol \lambda \in \mathbb R^{n+1}$ such that the following KKT conditions hold:
\begin{equation}\label{eq:KKT}
\begin{array}{rcll}
f(x_{i+1}) -f(x_{i-1}) + f'(x_i)(x_{i-1}-x_{i+1}) + \lambda_i - \lambda_{i-1} &= &0, & i = 1,\cdots,n,\\[5pt]
0 \, \leq \, x_{i+1} - x_i \, \perp \, \lambda_i & \geq & 0, & i = 0,\cdots,n.
\end{array}
\end{equation}
\begin{lemma}\label{lm:KKT_gradient}
Assume condition (a) holds. Let $(\boldsymbol x,\boldsymbol \lambda)$ be a pair of vectors satisfying (\ref{eq:KKT}), it holds that $\lambda_i = 0$ for all $i = 0,\cdots,n$.
\end{lemma}
\noindent\textbf{Proof.} Prove by contradiction. Suppose there exists an index $i^*$ such that $\lambda_{i*} > 0$. By the complementarity condition in (\ref{eq:KKT}), we have $x_{i^*} = x_{i^* + 1}$. \\
Claim (i): $a = x_0 = x_1 = \cdots = x_{i^*}.$ If $i^*=0$, Claim (i) holds readily. When $i^*\geq 1$,
by the equation in (\ref{eq:KKT}), we have
$$f(x_{i^*}) - f(x_{i^* -1}) + f'(x_{i^*})(x_{i^*-1} - x_{i^*}) + \lambda_{i^*} - \lambda_{i^*-1} \, = \, 0,$$ which implies that
$$\lambda_{i^*-1}\, = \, f(x_{i^*}) - f(x_{i^* -1}) + f'(x_{i^*})(x_{i^*-1} - x_{i^*}) + \lambda_{i^*}$$
By Mean Value Theorem, we have
$$\frac{f(x_{i^*}) - f(x_{i^* -1})}{x_{i^*} - x_{i^*-1}} \, = \, f'(\overline x),$$
for some $\overline x \in (x_{i^*-1},x_{i^*})$. By assumption (a), we have
$f'(\overline x) \geq f'(x_{i^*})$. Thus,
$$\frac{f(x_{i^*}) - f(x_{i^* -1})}{x_{i^*} - x_{i^*-1}} \geq f'(x_{i^*})$$ or equivalently,
$$f(x_{i^*}) - f(x_{i^* -1}) \geq f'(x_{i^*})(x_{i^*} - x_{i^*-1}). $$ Therefore,
$$\lambda_{i^*-1} \, = \, f(x_{i^*}) - f(x_{i^* -1}) + f'(x_{i^*})(x_{i^*-1} - x_{i^*}) + \lambda_{i^*} \, \geq \, \lambda_{i^*}\, > \, 0.$$
By applying the above arguments repeatedly, we can derive that $\lambda_i > 0$ for all $i\leq i^*$ and $i\geq 0$, and hence, $$a = x_0 = x_1 = \cdots = x_{i^*}.$$
Claim (ii): $x_{i^*} = x_{i^*+1} = \cdots = x_n =b.$ If $i^*=n-1$, then Claim (ii) holds readily. When $i^*\leq n-2$, from (\ref{eq:KKT}), we also have
$$f(x_{i^*+2}) - f(x_{i^*}) + f'(x_{i^*+1})(x_{i^*} - x_{i^*+2}) + \lambda_{i^*+1} - \lambda_{i^*} \, = \, 0,$$ which implies that
$$f(x_{i^*+2}) - f(x_{i^*+1}) + f'(x_{i^*+1})(x_{i^*+1} - x_{i^*+2}) + \lambda_{i^*+1} - \lambda_{i^*} \, = \, 0.$$
By Mean Value Theorem, we have
$$\frac{f(x_{i^*+2}) - f(x_{i^* +1})}{x_{i^*+2} - x_{i^*+1}} \, = \, f'(\tilde x),$$
for some $\tilde x \in (x_{i^*+1},x_{i^*+2})$. By assumption (a), we have
$f'(\tilde x) \leq f'(x_{i^*+1})$. Thus,
$$\frac{f(x_{i^*+2}) - f(x_{i^* +1})}{x_{i^*+2} - x_{i^*+1}} \, \leq \, f'(x_{i^*+1}),$$
or equivalently,
$$f(x_{i^*+2}) - f(x_{i^* +1}) + f'(x_{i^*+1}) ({x_{i^*+1} - x_{i^*+2}})  \leq \ 0.$$
Therefore,
$$\lambda_{i^*+1} \, = \, -\left[f(x_{i^*+2}) - f(x_{i^*+1}) + f'(x_{i^*+1})(x_{i^*+1} - x_{i^*+2}) \right] + \lambda_{i^*} \, \geq \, \lambda_{i^*} \, > \, 0.$$
If we keep applying this argument, we can obtain that $\lambda_i > 0$ for all $i\geq i^*$ and $i\leq n$, and hence, $$x_{i^*} = x_{i^*+1} = \cdots = x_n =b.$$ Since $a<b$, we have a contradiction.
This concludes the proof.\qed
\begin{corollary}\label{cor:no_binding}
Assume condition (a) holds. Let $(\boldsymbol x,\boldsymbol \lambda)$ be a pair of vectors satisfying (\ref{eq:KKT}), it holds that $x_0 \, < \, x_1 \, < \, \cdots\, <\, x_{n+1}$.
\end{corollary}
\noindent\textbf{Proof.} Prove by contradiction. Suppose there exists an index $i^*$ such that $x_{i^*} = x_{i^*+1}$ for some $i^*\in\{0,\cdots,n\}$. If $i^*=0$, by (\ref{eq:KKT}) and Lemma~\ref{lm:KKT_gradient}, we have
$$f(x_2)-f(x_0)+f'(x_1)(x_0-x_2) \, =\,0 .$$
Since $x_0=x_1$, we have
$$f(x_2)-f(x_1)+f'(x_1)(x_1-x_2) \, =\,0 .$$
Similar to the proof of Lemma~\ref{lm:KKT_gradient}, by applying Mean Value Theorem, we have
$$ f'(x_1) \, = \, \frac{f(x_2)-f(x_1)}{x_2-x_1} \, = \, f'(\widetilde x), $$
for some $\widetilde x \in (x_1,x_2)$. This contradicts with assumption (a). When $i^* \geq 1$, by (\ref{eq:KKT}) and Lemma~\ref{lm:KKT_gradient} we have
$$ f(x_{i^*})-f(x_{i^*-1}) = f'(x_{i^*})-f(x_{i^*-1}),$$
which also leads to a contradiction with assumption (a).\qed
\begin{corollary}\label{cor:gradient_0}
Assume condition (a) holds. Let $(\boldsymbol x,\boldsymbol \lambda)$ be a pair of vectors satisfying (\ref{eq:KKT}), it holds that $\nabla \phi(\boldsymbol x) = 0$.
\end{corollary}
\begin{corollary}\label{cor:quadratic}
If $f(x)=\alpha x^2 + \beta x + \gamma$ is a quadratic function with $\alpha < 0$, then the optimal solution of (\ref{eq:error_minimization}) is $$x_i = x_{i-1}+\frac{b-a}{n+1}, \,\, \forall \, i=1,\cdots,n.$$
\end{corollary}
\noindent\textbf{Proof.} By KKT conditions (\ref{eq:KKT}) and Lemma~\ref{lm:KKT_gradient}, we know that
$$f(x_{i+1}) -f(x_{i-1}) + f'(x_i)(x_{i-1}-x_{i+1}) \, =  \, 0, \,\, i = 1,\cdots,n.$$
Therefore
$$  2\alpha x_i + \beta  \, = \, f'(x_i) \, = \, \frac{f(x_{i+1}) -f(x_{i-1})}{x_{i+1}-x_{i-1}} \, = \, \frac{\alpha(x_{i+1}^2 - x_{i-1}^2) + \beta(x_{i+1} - x_{i-1})}{{x_{i+1}-x_{i-1}}}, \,\, i=1,\cdots,n.$$
This implies
$$x_i=\frac{x_{i-1}+ x_{i+1} }{2}, \,\,i=1,\cdots,n.$$ Hence, any $(\mathbf x,\mathbf \lambda)$ satisfying (\ref{eq:KKT}) must satisfy $x_1-x_0 = x_2-x_1 = \cdots = x_{n+1}-x_{n}$, which in turn implies that
$$x_i = x_{i-1}+\frac{b-a}{n+1}, \,\, \forall \, i=1,\cdots,n.$$
Notice that the optimization (\ref{eq:error_minimization}) has a continuous objective function and a convex closed bounded feasible region. Therefore, an optimal solution must exist and satisfy (\ref{eq:KKT}). Hence, the result holds readily.\qed
\begin{remark}
Corollary~\ref{cor:quadratic} states that when $f(x)$ is quadratic and concave, then the optimal distribution of the knots is locating them evenly in the interval of interests.
\end{remark}
\begin{lemma}\label{lm:no_ordering}
Assume condition (a) holds. Given a vector $\boldsymbol x \in \mathbb R^n$. If $\nabla \phi(\boldsymbol x) = 0$, and $a \leq x_i \leq b$, $i =1,\cdots,n$, then there exists $\boldsymbol \lambda \in \mathbb R^{n+1}$ so that $(\boldsymbol x,\boldsymbol \lambda)$ satisfies (\ref{eq:KKT}).
\end{lemma}
\noindent\textbf{Proof.} It suffices to show that $x_1 \leq x_2 \leq \cdots \leq x_n$. Then by letting $\boldsymbol \lambda = 0$, $(\boldsymbol x, \boldsymbol \lambda)$ satisfies (\ref{eq:KKT}) readily. Since $\nabla \phi(\boldsymbol x) =0$, we have
$$ f'(x_i) \, =\, \frac{f(x_{i+1})-f(x_{i-1})}{x_{i+1}-x_{i-1}}, \,\, \forall \, i=1,\cdots,n.$$
It is clear from the Mean Value Theorem that $$\frac{f(x_{i+1})-f(x_{i-1})}{x_{i+1}-x_{i-1}} = f'(\overline x),$$ for some $\overline x  \in (x_{i+1}, x_{i-1})$. By assumption (a), $f'(x)$ is a strictly decreasing function of $x$, therefore, we must have $x_{i-1} \leq x_i \leq x_{i+1}$ for all $i=1,\cdots,n$. This concludes the proof.
\qed

Combining Corollary~\ref{cor:gradient_0} and Lemma~\ref{lm:no_ordering}, we have the following theorem.
\begin{theorem}
Assume condition (a) holds. A pair of vectors $(\boldsymbol x, \boldsymbol \lambda)$ with $\boldsymbol x\in \mathbb R^n$ and $\boldsymbol \lambda \in \mathbb R^{n+1}$ satisfies KKT conditions~(\ref{eq:KKT}) if and only if $\nabla \phi(\boldsymbol x)=0$, $a \leq x_i \leq b$, and $\boldsymbol \lambda = 0$.
\end{theorem}

Next we look at sufficient conditions for optimality. The Hessian matrix of the objective function in (\ref{eq:error_minimization}) is given by
\begin{footnotesize}
\begin{equation}\label{eq:Hessian}
\left[
\begin{array}{ccccc}
(x_0-x_2)f''(x_1) & f'(x_2)-f'(x_1) & \cdots & 0 & 0 \\
f'(x_2)-f'(x_1) & (x_1-x_3)f''(x_2)& \cdots& 0 & 0\\
\cdots&\cdots&\cdots&\cdots& \cdots\\
0 & 0 & \cdots& (x_{n-2} - x_n)f''(x_{n-1}) & f'(x_n) - f'(x_{n-1}) \\
0 & 0 & \cdots& f'(x_n) - f'(x_{n-1})& (x_{n-1} - x_{n+1})f''(x_n)
\end{array}
\right]
\end{equation}
\end{footnotesize}
which is a tridiagonal matrix. This matrix is in general not positive definite. At a KKT point of (\ref{eq:error_minimization}), if the matrix (\ref{eq:Hessian}) is positive definite, then this KKT point must be a local maximum of (\ref{eq:error_minimization}). Applying a well known result regarding tridiagonal matrices, we have the following result whose proof is a simple application of a theorem from \cite{JNTsat96}.
\begin{proposition}
Let $x^*$ be a KKT point of (\ref{eq:error_minimization}), if it holds that for all $i=1,\cdots,n$ $$\left[f'(x^*_{i+1}) - f'(x^*_i)\right]^2 \, < \, \frac{1}{4} (x^*_{i-1}-x^*_{i+1}) (x^*_{i} - x^*_{i+2})f''(x^*_i)f''(x^*_{i+1})\frac{1}{\cos^2\left(\frac{\pi}{n+1}\right)},$$
then $x^*$ is a local minimum of (\ref{eq:error_minimization}).
\end{proposition}
\subsection{Approximating General Increasing Univariate Curves}\label{subsec:general}
We next move on to more general univariate curves. More specifically, we remove the concavity assumption. We assume there exists an $\epsilon>0$ such that
\begin{itemize}
\item[(b).] $f(x)$ is twice differentiable on $(a-\epsilon,b+\epsilon)$. Moreover, $f(x)>0$ and $f'(x)>0$ for all $x\in (a-\epsilon,b+\epsilon)$.
\end{itemize}
In this case, the PL approximation is not necessarily an under approximation and hence the error measurement in (\ref{eq:concave_measure}) is not appropriate anymore. In fact, (\ref{eq:concave_measure}) is the $\mathcal L_1$ norm of the difference between $f(x)$ and its PL approximation. Since the PL approximation is always an under approximation of $f(x)$ in the concave case, we are able to avoid taking the absolute value which leads to non-differentiability. In the general case, we therefore use the following error measure:
\begin{equation}\label{eq:error_general}
\sum_{i=0}^n \left[\int_{x_i}^{x_{i+1}} f(x) - \alpha_i x - \beta_i dx \right]^2
\end{equation}
For $i=0,\cdots,n$, we let
\begin{eqnarray}
\psi_i (x_i,x_{i+1}) & \triangleq & \left[\int_{x_i}^{x_{i+1}} f(x) - \alpha_i x - \beta_i dx \right]^2 \nonumber \\
& = & \left[ \int_{x_i}^{x_{i+1}} f(x)dx - \frac{1}{2}(f(x_i) + f(x_{i+1}))(x_{i+1} - x_i)\right]^2.\nonumber
\end{eqnarray}
The error minimization problem is then
\begin{equation}\label{eq:error_minimization_general}
\begin{array}{rl}
\displaystyle{\min_{x_1,\cdots,x_n}} & \displaystyle{ \sum_{i=1}^n \psi(x_i,x_{i+1})} \\[5pt]
\mbox{s.t.} & a \, = \, x_0 \, \leq \, x_1 \, \leq \, \cdots \, \leq \, x_n \, \leq \, x_{n+1} \, = \, b.
\end{array}
\end{equation}
Similar to the concave case, we introduce multipliers $\lambda_i$'s for the constraints.
\begin{equation}
\begin{array}{rlc}
\displaystyle{ \min_{x_1,\cdots, x_n} } & \displaystyle{ \sum_{i=0}^n \psi_i(x_i,x_{i+1}) }\\[10pt]
\mbox{s.t.} & x_0 - x_1 \, \leq \, 0 & (\lambda_0) \\
& x_1 - x_2 \, \leq \, 0 & (\lambda_1)\\
& \vdots & \vdots\\
& x_n - x_{n+1} \, \leq \, 0 & (\lambda_n)
\end{array}
\end{equation}
The Lagrangian function in this case is given by
\begin{equation}
L(\boldsymbol x,\boldsymbol \lambda) \, \triangleq \, \sum_{i=0}^n \psi_i(x_i,x_{i+1}) + \sum_{i=0}^n \lambda_i(x_i-x_{i+1}).
\end{equation}
Therefore, if $\boldsymbol x$ is an optimal solution of (\ref{eq:error_minimization_general}), then there exists $\boldsymbol \lambda \in \mathbb R^{n+1}$, such that:
\begin{equation}\label{eq:KKT_general}
\begin{array}{rcll}
\displaystyle { \frac{\partial \psi_{i-1}(x_{i-1},x_i)}{\partial x_i} + \frac{\partial \psi_{i}(x_i, x_{i+1})}{\partial x_i} } \, + \, \lambda_i \, - \, \lambda_{i-1} &= &0, & i = 1,\cdots,n,\\[5pt]
0 \, \leq \, x_{i+1} \, - \, x_i \, \perp \, \lambda_i & \geq & 0, & i = 0,\cdots,n.
\end{array}
\end{equation}
For each $i=0,\cdots,n$, the partial derivative of $\psi_i(x_i,x_{i+1})$ is given by
\begin{eqnarray}
\frac{\partial \psi_i(x_i,x_{i+1})}{\partial x_i} &=& -2 f(x_i) \int_{x_i}^{x_{i+1}} f(x) dx \, + \, f'(x_i)(x_{i+1}-x_i)\int_{x_i}^{x_{i+1}}f(x) dx \nonumber \\
&&- \left(f(x_i)+f(x_{i+1}) \right) \int_{x_i}^{x_{i+1}}f(x) dx \, - \,f(x_i) \left(f(x_i)+f(x_{i+1}) \right)(x_{i+1}-x_i)\nonumber\\
&&+ \frac{1}{2} \left(f(x_i)+f(x_{i+1}) \right) (x_{i+1}-x_i) \left[ f'(x_i)(x_{i+1}-x_i) -f(x_i)-f(x_{i+1})\right] \nonumber\\
& =& \left[f'(x_i)(x_{i+1}-x_i)-3f(x_i)-f(x_{i+1})\right]  \nonumber\\
&&\times \,\left[\int_{x_i}^{x_i+1} f(x) dx + \frac{1}{2}\left(f(x_i)+f(x_{i+1})\right)(x_{i+1}-x_i)\right]\label{eq:partial_x_i}
\end{eqnarray}
\begin{eqnarray}
\frac{\partial \psi_i(x_i,x_{i+1})}{\partial x_{i+1}} &=& 2 f(x_{i+1}) \int_{x_i}^{x_{i+1}} f(x) dx  \,+ \, f'(x_{i+1})(x_{i+1}-x_i)\int_{x_i}^{x_{i+1}}f(x) dx \nonumber \\
  && \nonumber\\
  && + \left(f(x_i)+f(x_{i+1}) \right) \int_{x_i}^{x_{i+1}}f(x) dx \, + \, f(x_{i+1}) \left(f(x_i)+f(x_{i+1}) \right)(x_{i+1}-x_i) \nonumber \\
  && \frac{1}{2} \left(f(x_i)+f(x_{i+1}) \right) (x_{i+1}-x_i) \left[ f'(x_{i+1})(x_{i+1}-x_i) +f(x_i)+f(x_{i+1})\right] \nonumber\\
  & =& \left[f'(x_i)(x_{i+1}-x_i)+f(x_i)+3f(x_{i+1})\right]  \nonumber\\
&&\times \,\left[\int_{x_i}^{x_i+1} f(x) dx + \frac{1}{2}\left(f(x_i)+f(x_{i+1})\right)(x_{i+1}-x_i)\right]\label{eq:partial_x_i+1}
\end{eqnarray}
Under assumption (b) it is easy to see that $$\frac{\partial \psi_i(x_i,x_{i+1})}{\partial x_{i+1}} \, \geq \, 0$$ and $$\frac{\partial \psi_i(x_i,x_{i+1})}{\partial x_{i+1}} \, =\, 0 \,\Longleftrightarrow \, x_i \, = \,x_{i+1}.$$
On the other hand,
$$x_i \, = \,x_{i+1} \, \Longrightarrow \, \frac{\partial \psi_i(x_i,x_{i+1})}{\partial x_{i}} \, =\, 0 .$$
\begin{lemma}\label{lm:general}
Assume condition (b) holds. Let $(\boldsymbol x,\boldsymbol \lambda)$ be a pair of vectors satisfying (\ref{eq:KKT}), it holds that $x_n<x_{n+1}$.
\end{lemma}
\noindent\textbf{Proof.} Prove by contradiction. If $x_n=x_{n+1}$, by (\ref{eq:KKT_general}) we have
$$\lambda_{n-1} \, = \, \lambda_{n} + \frac{\partial \psi_{n-1}(x_{n-1},x_n)}{\partial x_n} + \frac{\partial \psi_{n}(x_n, x_{n+1})}{\partial x_n} \, = \, \lambda_{n} + \frac{\partial \psi_{n-1}(x_{n-1},x_n)}{\partial x_n} \, \geq \, \frac{\partial \psi_{n-1}(x_{n-1},x_n)}{\partial x_n}.$$
Now, if $\lambda_{n-1} > 0$, then $x_{n-1}=x_n$. Else if $\lambda_{n-1} = 0$, then we must have $$\frac{\partial \psi_{n-1}(x_{n-1},x_n)}{\partial x_n} \, \Longrightarrow \, x_n=x_{n-1}.$$ Therefore, we must also have $x_n = x_{n-1}$. If we keep applying this argument, we derive that $$b = x_{n+1} = \cdots =x_0 = a,$$ which is a contradiction.\qed
\begin{remark}
Lemma \ref{lm:general} ascertains that the last knot does not coincide with the end point of the interval of interests in any local optimal solutions. Therefore, the transformation we discuss in the next section works for this case also. \qed
\end{remark}
\section{A Gradient Projection Algorithm}\label{sec:alg}
To solve the optimization problem (\ref{eq:error_minimization}), we propose to apply gradient projection method due to its simplicity. Moreover, as we see below, we can modify our formulation so that the feasible set allows a strongly polynomial time projection algorithm. We notice that the feasible set of (\ref{eq:error_minimization}) is a nonnegative monotone cone restricted by an upper bound. On the other hand, there is a strongly polynomial time algorithm to project onto a monotone nonnegative cone \cite{nemeth2012project}. Therefore, we rewrite the optimization problem~(\ref{eq:error_minimization}) so that the feasible region is a monotone nonnegative cone. Note that we use (\ref{eq:error_minimization}) to demonstrate our approach here. Lemma~\ref{lm:general} ensures that the same approach can be applied to (\ref{eq:error_minimization_general}). We let
\begin{equation}\label{eq:y_transformation}
y_i \, = \, \frac{x_i-a}{b-x_i}.
\end{equation}
We notice that $y_i$'s are well defined when $x_i \neq b$. By Corollary~\ref{cor:no_binding}, for all KKT point $(\boldsymbol x)$ of (\ref{eq:error_minimization}), $x_i<b$ for all $i=1,\cdots,n$. Therefore, after introducing $y_i$'s, we can rewrite (\ref{eq:error_minimization}) in terms of $y_i$'s without changing optimal solutions of (\ref{eq:error_minimization}). In fact, from equation~(\ref{eq:y_transformation}) we have
$$x_i \, = \, \frac{by_i +a}{1+y_i} \, = \, b - \frac{b-a}{1+y_i}.$$
Let $y_0=0$ and $\boldsymbol y = [y_1,\cdots,y_n]^T$, the objective function in (\ref{eq:error_minimization}) becomes
\begin{eqnarray}
\Phi (\boldsymbol y) & \triangleq & -\frac{1}{2}\sum_{i=0}^{n-1} \left(\frac{b-a}{1+y_{i}}-\frac{b -a}{1+y_{i+1}}\right) \left(f\left(\frac{b y_{i+1}+a}{1+y_{i+1}}\right)+ f\left(\frac{b y_i+a}{1+y_i}\right)\right) \nonumber\\
&& \,- \,\frac{1}{2} \left(b - \frac{b y_n+a}{1+y_n}\right)\left( f(b) + f\left(\frac{b y_n+a}{y_n+1}\right)\right)\nonumber
\end{eqnarray}
The optimization problem~(\ref{eq:error_minimization}) can be rewritten as
\begin{equation}\label{eq:error_minimization_y}
\begin{array}{rl}
\displaystyle{ \min_{y_1,\cdots, y_n} } & \Phi(\boldsymbol y)\\[5pt]
\mbox{s.t.} & 0\, \leq y_1 \, \leq \, \cdots\, \leq \, y_n.
\end{array}
\end{equation}
Let the monotone nonnegative cone in $\mathbb R^n$ be denoted by $\mathcal M^n$. For any vector $\boldsymbol x \in \mathbb R^n$, the projection on to $\mathcal M^n$, denoted by $\Pi_{\mathcal M^n}(\boldsymbol x)$, is the optimal solution of the following quadratic program:
\begin{equation}\label{eq:projection}
\begin{array}{rcrl}
\Pi_{\mathcal M^n}(\boldsymbol x) & = & \displaystyle{\min_{\boldsymbol y \in \mathbb R^n} }&(\boldsymbol x-\boldsymbol y)^T(\boldsymbol x-\boldsymbol y) \\[5pt]
&&\mbox{s.t.} & 0\leq y_1 \leq \cdots \leq y_n.
\end{array}
\end{equation}
As shown in \cite{nemeth2012project}, for any $\boldsymbol x\in \mathcal R^n$, $\Pi_{\mathcal M^n}(\boldsymbol x)$ is computable in strongly polynomial time. Therefore, we apply the Spectral Projected Gradient (SPG) Algorithm \cite{birgin2000nonmonotone} to solve optimization problem~(\ref{eq:error_minimization_y}).
\begin{algorithm}[htbp!]
	\caption{Spectral Projected Gradient Algorithm}\label{alg:SPG}
	\begin{algorithmic}[1]
		\State Given $ \boldsymbol y_0 $, step bounds $ 0 < \alpha_{\text{min}} < \alpha_{\text{max}} $
		\State Initial step length $\alpha_{bb} \in [\alpha_{\text{min}}, \alpha_{\text{max}}]$, and history length $ h $
		\State \textbf{While} not converge:
		\State \indent $\bar{\alpha}_k = \min\{ \alpha_{\text{max}}, \max\{ \alpha_{\text{min}}, \alpha_{bb}\} \} $
		\State \indent $d_k = \Pi_{\mathcal M^n}(\boldsymbol y_k - \bar{\alpha}_k \nabla \Phi(\boldsymbol y_k)) - \boldsymbol y_k $.
		\State \indent Set bound $f_b = \max\{ \Phi(\boldsymbol y_k), \cdots, \Phi(\boldsymbol y_{k - h}) \}$
		\State \indent $ \alpha = 1  $
		\State \indent \textbf{While} $ \Phi(\boldsymbol y_k + \alpha d_k) > f_b + \nu \alpha \nabla \Phi (\boldsymbol y_k)^{T} d_k$:
		\State \indent \indent Select $ \alpha $ randomly from Uniform distribution $ U(0, \alpha) $.
		\State \indent $ \boldsymbol y_{k+1} = \boldsymbol y_k + \alpha d_k $
		\State \indent $ s_{k} = \boldsymbol y_{k+1} - \boldsymbol y_k $
		\State \indent $ z_k = \nabla \Phi(\boldsymbol y_{k+1}) - \nabla \Phi(\boldsymbol y_k) $
		\State \indent $ \alpha_{bb} = \boldsymbol z_k^{T}\boldsymbol z_k/s_k^{T}s_k $
		\State \indent $ k = k + 1  $	
	\end{algorithmic}
\end{algorithm}


The termination conditions we use in our experiments include:
\begin{itemize}
\item[(1)] maximum number of iterations reached; or
\item[(2)] not enough improvement; or
\item[(3)] optimality condition is satisfied, i.e., $\|d_k\|\leq \epsilon$ for some small predefined $\epsilon>0$.
\end{itemize}
\section{Numerical Experiments}\label{sec:exp}
In our numerical experiments, we consider 5 different types of curves as follows:
\begin{itemize}
\item[Type 1:] Generalized Logtistic:
\begin{equation}
f_1(x) \,  = \, v_1 + \frac{v_2}{(1+s e^{d_1x + d_2})^{\frac{1}{s}}}.
\end{equation}
\item[Type 2:] Gompertz:
\begin{equation}
f_2(x) \, = \, v_1 + v_2 e^{s e^{d_1x+d_2}}.
\end{equation}
\item[Type 3:] Weibull:
\begin{equation}
f_3(x) \, = \, v_1 + v_2 e^{-(d_1 x + d_2)^s}.
\end{equation}
\item[Type 4:] Arctangent:
\begin{equation}
f_4(x) \, = \, v_1 + v_2 \arctan(d_1 x+ d_2).
\end{equation}
\item[Type 5:] Algebraic:
\begin{equation}
f_5(x) \, = \, v_1 + \frac{v_2}{(d_1 x^s + d_2)^{\frac{1}{s}}}.
\end{equation}
\end{itemize}
All of the above 5 types of curves can be either C-shaped or non-C-shaped depending on the values of the parameters. In the following table, we list all the curves we used in our experiments with the parameter values and intervals of interest.
\begin{table}[htbp!]
\centering
\begin{tabular}{|c|c|c|c|c|c|c|c|c|c|}
  \hline
  Name & Type & $v_1$ & $v_2$ & $s$ & $d_1$ & $d_2$ & Concave & $a$ & $b$ \\
  \hline
logistic1a & Logistic & 0.0 & 1.0 & 1.0 & -1.0 & 0.0 & Y & 0.0 & 2.0\\
logistic2a & Logistic & 0.0 & 1.2 & 1.0 & -1.5 & 0.0 & Y & 0.0 & 2.0\\
logistic3a & Logistic & 0.0 & 1.2 & 1.0 & -1.5 & 0.0 & Y & 0.25 & 1.75\\
gompertz1a & Gompertz & 0.0 & 1.0 & -1.0 & -1.0 &0.0 & Y & 0.0 & 6.0\\
weibull1a & Weibull & 1.0 & -1.0 & 2.0 & 1.0 & 0.0 & Y & -2.0 & -0.5\\
weibull2a & Weibull & 1.0 & -1.0 & 4.0 & 1.0 & 0.0 & Y & 1.0 & 3.0\\
weibull3a & Weibull & 1.0 & -1.0 & 2.2 & 1.0 & 0.0 & Y & 1.0 & 3.0\\
logistic1b & Logistic & 0.0 & 1.0 & 1.0 & -1.0 & 0.0 & N & -2.0 & 2.0\\
logistic2b & Logistic & 0.0 & 1.2 & 1.0 & -1.5 & 0.0 & N & -2.0 & 2.0\\
logistic3b & Logistic & 0.0 & 1.2 & 1.0 & -1.5 & 0.0 & N & -1.75 & 1.75\\
gompertz1b & Gompertz & 0.0 & 1.0 & -1.0 & -1.0 & 0.0 & N & -3.0 & 6.0\\
gompertz2b & Gompertz & 0.0 & 2.0 & -0.5 & 0.5 & 0.0 & N & -6.0 & 6.0\\
gompertz3b & Gompertz & 0.0 & 2.0 & -3.0 & 0.5 & 0.0 & N & -6.0 & 6.0\\
weibull1b & Weibull & 1.0 & -1.0 & 2.0 & 1.0 & 0.0 & N & -2.0 & 2.0\\
weibull2b & Weibull & 1.0 & -1.0 & 4.0 & 1.0 & 0.0 & N & -2.0 & 2.0\\
arctan1b & Arctan & 0.0 & 1.0 & - & 1.0 & 0.0 & N & -6.0 & 6.0\\
arctan2b & Arctan & 0.0 & 1.0 & - & 0.5 & 0.5 & N & -6.0 & 6.0\\
arctan3b & Arctan & 0.0 & -4.0 & - & 0.5 & 0.0 & N & -6.0 & 6.0\\
algebraic1b & Algebraic & 1.0 & 2.0 & 2.0 & 1.5 & 2.0 & N & -2.0 & 4.0\\
algebraic2b & Algebraic & 2.0 & 2.0 & 2.0 & 1.5 & 2.0 & N & -4.0 & 2.0\\
  \hline
\end{tabular}
\caption{List of curves in numerical experiments, $a$ and $b$ are the lower and upper bounds of the interval of interests.} \label{tab:list_curve}
\end{table}
As we can see in Table \ref{tab:list_curve}, the first 7 curves are concave, and hence we apply reformulation (\ref{eq:error_minimization_y}) of (\ref{eq:error_minimization}) to locate the knots. We apply a sequential quadratic programming (SQP) solver implemented in Scipy optimize module (SciPy.optimize.minimze(method='SLSQP')). See \cite{Kraft88} for details of this algorithm, and the SPG algorithm presented in the previous section to solve the formulation. For each curve, we conduct experiments with 4 and 8 knots, respectively. For both of the algorithms we initialize them with equally distributed knots. Our results are presented in Table \ref{tab:concave curves}, the initial approximation error measured by (\ref{eq:concave_measure}) (with equally distributed knots) is given in column ``orig error". The approximation errors, measured by (\ref{eq:concave_measure}), of the solutions found by our SPG algorithm and SciPy SQP algorithm are given in columns ``SPG error" and ``SQP error", respectively. To compare the results of the two algorithm, we also include the relative difference in the error of the two algorithm in column ``diff in error". More specifically, we compute
$$ \mbox{diff in error} \, =\, \frac{\mbox{SQP error} - \mbox{SPG error}}{\mbox{SQP error}} \times 100 \%. $$
As we can see from Table \ref{tab:concave curves}, both the SPG algorithm and the SQP algorithm successfully reduce the approximation error. This shows the effectiveness of our overall approach. We also observe that when the number of knots is small, the reduction is more significant. This observation is consistent with our intuition. We notice that in most of the cases (8 out of 14), the error of the results produced by SPG is smaller than the results produced by SQP. We also observed that SPG algorithm typically is slower than the SQP algorithm. One of the reasons is that the SPG algorithm is completely written in Python and is home grown and not optimized while the SQP implementation in SciPy is essentially a wrapper of a piece of well-developed code written in Fortran. Another reason for this performance difference is that the SQP algorithm uses second order approximation information (e.g. quasi-Newton's method) while the SPG algorithm uses only first order information.
\begin{table}[htbp!]
  \centering
  \begin{tabular}{|c|c|c|c|c|c|c|c|}
  \hline
  curve name & $a$	& $b$ & \# knots & orig error & SPG error & SQP error & diff in error \% \\
  \hline
  logistic1a & 0.00 & 2.00 & 4 & 6.166057E-07 & 2.925162E-08 & 3.811412E-07 & 9.2325E+01\\
  logistic1a & 0.00 & 2.00 & 8 & 3.901868E-08 & 3.061227E-08 & 3.901868E-08 & 2.1545E+01\\
  logistic2a & 0.00 & 2.00 & 4 & 4.546293E-06 & 2.354395E-07 & 2.872030E-06 & 9.1802E+01\\
  logistic2a & 0.00 & 2.00 & 8 & 2.704366E-07 & 2.204636E-07 & 1.926915E-07 & -1.4413E+01\\
  logistic3a & 0.25 & 1.75 & 4 & 8.866112E-07 & 4.936710E-08 & 6.979750E-07 & 9.2927E+01\\
  logistic3a & 0.25 & 1.75 & 8 & 5.594112E-08 & 5.129870E-08 & 5.594112E-08 & 8.2988E+00\\
  gompertz1a & 0.00 & 6.00 & 4 & 3.319009E-04 & 1.414356E-05 & 1.235959E-04 & 8.8557E+01\\
  gompertz1a & 0.00 & 6.00 & 8 & 3.075644E-05 & 2.042979E-05 & 7.922744E-06 & -1.5786E+02\\
  weibull1a & -2.00 & -0.50 & 4 & 8.351922E-06 & 3.460830E-07 & 3.036610E-06 & 8.8603E+01\\
  weibull1a & -2.00 & -0.50 & 8 & 4.678674E-07 & 3.435558E-07 & 1.621089E-07 & -1.1193E+02\\
  weibull2a & 1.00 & 3.00 & 4 & 6.853906E-06 & 4.671216E-06 & 1.566036E-06 & -1.9828E+02\\
  weibull2a & 1.00 & 3.00 & 8 & 7.173659E-06 & 2.799612E-06 & 5.447805E-07 & -4.1390E+02\\
  weibull3a & 1.00 & 3.00 & 4 & 1.647924E-05 & 1.160405E-06 & 4.781666E-06 & 7.5732E+01\\
  weibull3a & 1.00 & 3.00 & 8 & 1.654462E-06 & 1.462703E-06 & 2.952623E-07 & -3.9539E+02\\
  \hline
  \end{tabular}
  \caption{Numerical results on concave curves}\label{tab:concave curves}
\end{table}
For all the curves, including the concave ones and the non-concave ones, we apply the same two algorithms on the reformulation (\ref{eq:error_minimization_y}) of (\ref{eq:error_minimization_general}). For each curve, we also conduct experiments with 4 and 8 knots respectively. We initialize both algorithms with equally distributed knots. Our results are presented in Table \ref{tab: non_concave_results}. The initial approximation error measured by (\ref{eq:error_general}) (with equally distributed knots) is given in column ``orig error". The approximation errors, measured by (\ref{eq:error_general}), of the solutions found by our SPG algorithm and SciPy SQP algorithm are given in columns ``SPG error" and ``SQP error", respectively. We also include the difference in error in column ``diff in error". Our observations are very similar to the concave case. The results again demonstrate the effectiveness of the proposed approach in the non-concave case.
\begin{table}[htbp!]
  \centering
  \begin{tabular}{|c|c|c|c|c|c|c|c|}
    \hline
    curve name & $a$ & b & \# knots & orig error &	SPG error & SQP error & diff in error \% \\
    \hline
    logistic1a & 0.00 & 2.00 & 4 & 6.166057E-07 & 3.901868E-08 & 6.166057E-07 & 9.3672E+01\\
    logistic1a & 0.00 & 2.00 & 8 & 3.901868E-08 & 3.901868E-08 & 3.901868E-08 & -2.8153E-12\\
    logistic2a & 0.00 & 2.00 & 4 & 4.546293E-06 & 2.702524E-07 & 4.546293E-06 & 9.4056E+01\\
    logistic2a & 0.00 & 2.00 & 8 & 2.704366E-07 & 2.702524E-07 & 2.704366E-07 & 6.8089E-02\\
    logistic3a & 0.25 & 1.75 & 4 & 8.866112E-07 & 5.594112E-08 & 8.866112E-07 & 9.3690E+01\\
    logistic3a & 0.25 & 1.75 & 8 & 5.594112E-08 & 5.594112E-08 & 5.594112E-08 & 0.0000E+00\\
    gompertz1a & 0.00 & 6.00 & 4 & 3.319009E-04 & 2.113859E-05 & 3.199724E-04 & 9.3394E+01\\
    gompertz1a & 0.00 & 6.00 & 8 & 3.075644E-05 & 2.293567E-05 & 3.075644E-05 & 2.5428E+01\\
    weibull1a & -2.00 & -0.50 & 4 & 8.351922E-06 & 4.385014E-07 & 8.351922E-06 & 9.4750E+01\\
    weibull1a & -2.00 & -0.50 & 8 & 4.678674E-07 & 4.385014E-07 & 4.678674E-07 & 6.2766E+00\\
    weibull2a & 1.00 & 3.00 & 4 & 6.853906E-06 & 1.219792E-06 & 6.076741E-05 & 9.7993E+01\\
    weibull2a & 1.00 & 3.00 & 8 & 7.173659E-06 & 1.219815E-06 & 7.173659E-06 & 8.2996E+01\\
    weibull3a & 1.00 & 3.00 & 4 & 1.647924E-05 & 1.049252E-06 & 1.647924E-05 & 9.3633E+01\\
    weibull3a & 1.00 & 3.00 & 8 & 1.654462E-06 & 1.049256E-06 & 1.654462E-06 & 3.6580E+01\\
    logistic1b & -2.00 & 2.00 & 4 & 2.287906E-05 & 6.588572E-07 & 2.287906E-05 & 9.7120E+01\\
    logistic1b & -2.00 & 2.00 & 8 & 2.049227E-06 & 6.588572E-07 & 2.049227E-06 & 6.7848E+01\\
    logistic2b & -2.00 & 2.00 & 4 & 2.232474E-04 & 1.871111E-06 & 4.384989E-10 & -4.2661E+05\\
    logistic2b & -2.00 & 2.00 & 8 & 1.593240E-05 & 1.855617E-06 & 1.593240E-05 & 8.8353E+01\\
    logistic3b & -1.75 & 1.75 & 4 & 9.481086E-05 & 6.294456E-07 & 9.481086E-05 & 9.9336E+01\\
    logistic3b & -1.75 & 1.75 & 8 & 7.285415E-06 & 6.294240E-07 & 7.285415E-06 & 9.1360E+01\\
    gompertz1b & -3.00 & 6.00 & 4 & 7.738086E-03 & 2.234425E-04 & 1.764598E-05 & -1.1663E+03\\
    gompertz1b & -3.00 & 6.00 & 8 & 7.514605E-04 & 4.389726E-04 & 3.334459E-04 & -3.1647E+01\\
    gompertz2b & -6.00 & 6.00 & 4 & 2.285238E-02 & 4.149349E-04 & 1.900082E-06 & -2.1738E+04\\
    gompertz2b & -6.00 & 6.00 & 8 & 1.720082E-03 & 3.717897E-04 & 6.304516E-04 & 4.1028E+01\\
    gompertz3b & -6.00 & 6.00 & 4 & 2.352946E-02 & 2.658223E-04 & 1.200728E-05 & -2.1138E+03\\
    gompertz3b & -6.00 & 6.00 & 8 & 1.473251E-03 & 3.313604E-04 & 1.510322E-04 & -1.1940E+02\\
    weibull1b & -2.00 & 2.00 & 4 & 6.166059E-03 & 1.046891E-04 & 7.516882E-08 & -1.3917E+05\\
    weibull1b & -2.00 & 2.00 & 8 & 4.069463E-04 & 1.461744E-04 & 2.563262E-05 & -4.7027E+02\\
    weibull2b & -2.00 & 2.00 & 4 & 6.091507E-03 & 6.342666E-04 & 6.682868E-07 & -9.4809E+04\\
    weibull2b & -2.00 & 2.00 & 8 & 1.316705E-03 & 5.488979E-04 & 3.087013E-05 & -1.6781E+03\\
    arctan1b & -6.00 & 6.00 & 4 & 4.205023E-02 & 3.012734E-03 & 2.604656E-13 & -1.1567E+12\\
    arctan1b & -6.00 & 6.00 & 8 & 1.080821E-02 & 1.123970E-03 & 3.535369E-04 & -2.1792E+02\\
    arctan2b & -6.00 & 6.00 & 4 & 5.327812E-02 & 5.703777E-04 & 7.732105E-06 & -7.2767E+03\\
    arctan2b & -6.00 & 6.00 & 8 & 2.619283E-03 & 6.882992E-04 & 1.185626E-04 & -4.8054E+02\\
    arctan3b & -6.00 & 6.00 & 4 & 4.515495E-01 & 9.070453E-03 & 7.349335E-07 & -1.2341E+06\\
    arctan3b & -6.00 & 6.00 & 8 & 4.121905E-02 & 1.109400E-02 & 2.393317E-05 & -4.6254E+04\\
    algebraic1b & -2.00 & 4.00 & 4 & 9.546650E-02 & 2.381183E-03 & 1.824035E-02 & 8.6946E+01\\
    algebraic1b & -2.00 & 4.00 & 8 & 5.375949E-03 & 3.003785E-03 & 1.225231E-03 & -1.4516E+02\\
    algebraic2b & -4.00 & 2.00 & 4 & 9.546650E-02 & 2.239908E-03 & 1.824034E-02 & 8.7720E+01\\
    algebraic2b & -4.00 & 2.00 & 8 & 5.375949E-03 & 1.769499E-03 & 1.225233E-03 & -4.4421E+01\\
    \hline
  \end{tabular}
  \caption{Numerical results on non-concave curves.}\label{tab: non_concave_results}
\end{table}
\section{Conclusion}\label{sec:con}
In this paper, we studied the piecewise linear approximation of univariate nonlinear function via the optimal location of knots. Given the number of knots, we formulate optimization problems to find the optimal knot location so that the PL approximation error is minimized. Properties of the optimization problems were studied, and reformulations of the original problems were derived based on their properties. The reformulations allowed us to apply a simple projection algorithm to solve the optimization problem. We demonstrated the efficiency of the proposed approach with extensive numerical experiments. Possible future research includes developing a convex measure of the approximation error so that global optimality can be achieved. Another possible area of research is the development of more sophisticated algorithms to solve the knot locating problem.

\end{document}